
\magnification=1200

\hsize=14cm
\vsize=20.5cm
\parindent=0cm   \parskip=0pt
\pageno=1

\def\ind{\hskip 1cm\relax}

\ifnum\mag=\magstep1
\hoffset=-0.5cm   
\voffset=-0.5cm   
\fi


\pretolerance=500 \tolerance=1000  \brokenpenalty=5000

\catcode`\@=11

\font\eightrm=cmr8         \font\eighti=cmmi8
\font\eightsy=cmsy8        \font\eightbf=cmbx8
\font\eighttt=cmtt8        \font\eightit=cmti8
\font\eightsl=cmsl8        \font\sixrm=cmr6
\font\sixi=cmmi6           \font\sixsy=cmsy6
\font\sixbf=cmbx6


\font\tengoth=eufm10       \font\tenbboard=msbm10
\font\eightgoth=eufm8      \font\eightbboard=msbm8
\font\sevengoth=eufm7      \font\sevenbboard=msbm7
\font\sixgoth=eufm6        \font\fivegoth=eufm5

\font\tencyr=wncyr10       
\font\eightcyr=wncyr8      
\font\sevencyr=wncyr7      
\font\sixcyr=wncyr6


\skewchar\eighti='177 \skewchar\sixi='177
\skewchar\eightsy='60 \skewchar\sixsy='60


\newfam\gothfam           \newfam\bboardfam
\newfam\cyrfam

\def\tenpoint{%
  \textfont0=\tenrm \scriptfont0=\sevenrm \scriptscriptfont0=\fiverm
  \def\rm{\fam\z@\tenrm}%
  \textfont1=\teni  \scriptfont1=\seveni  \scriptscriptfont1=\fivei
  \def\oldstyle{\fam\@ne\teni}\let\old=\oldstyle
  \textfont2=\tensy \scriptfont2=\sevensy \scriptscriptfont2=\fivesy
  \textfont\gothfam=\tengoth \scriptfont\gothfam=\sevengoth
  \scriptscriptfont\gothfam=\fivegoth
  \def\goth{\fam\gothfam\tengoth}%
  \textfont\bboardfam=\tenbboard \scriptfont\bboardfam=\sevenbboard
  \scriptscriptfont\bboardfam=\sevenbboard
  \def\bb{\fam\bboardfam\tenbboard}%
 \textfont\cyrfam=\tencyr \scriptfont\cyrfam=\sevencyr
  \scriptscriptfont\cyrfam=\sixcyr
  \def\cyr{\fam\cyrfam\tencyr}%
  \textfont\itfam=\tenit
  \def\it{\fam\itfam\tenit}%
  \textfont\slfam=\tensl
  \def\sl{\fam\slfam\tensl}%
  \textfont\bffam=\tenbf \scriptfont\bffam=\sevenbf
  \scriptscriptfont\bffam=\fivebf
  \def\bf{\fam\bffam\tenbf}%
  \textfont\ttfam=\tentt
  \def\tt{\fam\ttfam\tentt}%
  \abovedisplayskip=12pt plus 3pt minus 9pt
  \belowdisplayskip=\abovedisplayskip
  \abovedisplayshortskip=0pt plus 3pt
  \belowdisplayshortskip=4pt plus 3pt
  \smallskipamount=3pt plus 1pt minus 1pt
  \medskipamount=6pt plus 2pt minus 2pt
  \bigskipamount=12pt plus 4pt minus 4pt
  \normalbaselineskip=12pt
  \setbox\strutbox=\hbox{\vrule height8.5pt depth3.5pt width0pt}%
  \let\bigf@nt=\tenrm       \let\smallf@nt=\sevenrm
  \normalbaselines\rm}

\def\eightpoint{%
  \textfont0=\eightrm \scriptfont0=\sixrm \scriptscriptfont0=\fiverm
  \def\rm{\fam\z@\eightrm}%
  \textfont1=\eighti  \scriptfont1=\sixi  \scriptscriptfont1=\fivei
  \def\oldstyle{\fam\@ne\eighti}\let\old=\oldstyle
  \textfont2=\eightsy \scriptfont2=\sixsy \scriptscriptfont2=\fivesy
  \textfont\gothfam=\eightgoth \scriptfont\gothfam=\sixgoth
  \scriptscriptfont\gothfam=\fivegoth
  \def\goth{\fam\gothfam\eightgoth}%
  \textfont\cyrfam=\eightcyr \scriptfont\cyrfam=\sixcyr
  \scriptscriptfont\cyrfam=\sixcyr
  \def\cyr{\fam\cyrfam\eightcyr}%
  \textfont\bboardfam=\eightbboard \scriptfont\bboardfam=\sevenbboard
  \scriptscriptfont\bboardfam=\sevenbboard
  \def\bb{\fam\bboardfam}%
  \textfont\itfam=\eightit
  \def\it{\fam\itfam\eightit}%
  \textfont\slfam=\eightsl
  \def\sl{\fam\slfam\eightsl}%
  \textfont\bffam=\eightbf \scriptfont\bffam=\sixbf
  \scriptscriptfont\bffam=\fivebf
  \def\bf{\fam\bffam\eightbf}%
  \textfont\ttfam=\eighttt
  \def\tt{\fam\ttfam\eighttt}%
  \abovedisplayskip=9pt plus 3pt minus 9pt
  \belowdisplayskip=\abovedisplayskip
  \abovedisplayshortskip=0pt plus 3pt
  \belowdisplayshortskip=3pt plus 3pt
  \smallskipamount=2pt plus 1pt minus 1pt
  \medskipamount=4pt plus 2pt minus 1pt
  \bigskipamount=9pt plus 3pt minus 3pt
  \normalbaselineskip=9pt
  \setbox\strutbox=\hbox{\vrule height7pt depth2pt width0pt}%
  \let\bigf@nt=\eightrm     \let\smallf@nt=\sixrm
  \normalbaselines\rm}

\tenpoint


\def\pc#1{\bigf@nt#1\smallf@nt}         \def\pd#1 {{\pc#1} }


\catcode`\;=\active
\def;{\relax\ifhmode\ifdim\lastskip>\z@\unskip\fi
\kern\fontdimen2  -1.2 \fontdimen3 \string;}

\catcode`\:=\active
\def:{\relax\ifhmode\ifdim\lastskip>\z@\unskip\fi\penalty\@M\ \fi\string:}

\catcode`\!=\active
\def!{\relax\ifhmode\ifdim\lastskip>\z@
\unskip\fi\kern\fontdimen2  -1.1 \fontdimen3 \string!}

\catcode`\?=\active
\def?{\relax\ifhmode\ifdim\lastskip>\z@
\unskip\fi\kern\fontdimen2  -1.1 \fontdimen3 \string?}

\def\^#1{\if#1i{\accent"5E\i}\else{\accent"5E #1}\fi}
\def\"#1{\if#1i{\accent"7F\i}\else{\accent"7F #1}\fi}

\frenchspacing


\newtoks\auteurcourant      \auteurcourant={\hfil}
\newtoks\titrecourant       \titrecourant={\hfil}

\newtoks\hautpagetitre      \hautpagetitre={\hfil}
\newtoks\baspagetitre       \baspagetitre={\hfil}

\newtoks\hautpagegauche
\hautpagegauche={\eightpoint\rlap{\folio}\hfil\the\auteurcourant\hfil}
\newtoks\hautpagedroite
\hautpagedroite={\eightpoint\hfil\the\titrecourant\hfil\llap{\folio}}

\newtoks\baspagegauche      \baspagegauche={\hfil}
\newtoks\baspagedroite      \baspagedroite={\hfil}

\newif\ifpagetitre          \pagetitretrue

\def\raggedbottom{\topskip 10pt plus 36pt\r@ggedbottomtrue}

\long\def\ctitre#1\endctitre{%
    \ifdim\lastskip<24pt\vskip-\lastskip\bigbreak\bigbreak\fi
  		\vbox{\parindent=0pt\leftskip=0pt plus 1fill
          \rightskip=\leftskip
          \parfillskip=0pt\bf#1\par}
    \bigskip\nobreak}

\long\def\section#1\endsection{%
\vskip 0pt plus 3\normalbaselineskip
\penalty-250
\vskip 0pt plus -3\normalbaselineskip
\Bigbreak
\message{[section \string: #1]}{\bf#1\unskip}\pointir}

\def\signature#1\endsignature{\vskip 15mm minus 5mm\rightline{\vtop{#1}}}

\mathcode`A="7041 \mathcode`B="7042 \mathcode`C="7043 \mathcode`D="7044
\mathcode`E="7045 \mathcode`F="7046 \mathcode`G="7047 \mathcode`H="7048
\mathcode`I="7049 \mathcode`J="704A \mathcode`K="704B \mathcode`L="704C
\mathcode`M="704D \mathcode`N="704E \mathcode`O="704F \mathcode`P="7050
\mathcode`Q="7051 \mathcode`R="7052 \mathcode`S="7053 \mathcode`T="7054
\mathcode`U="7055 \mathcode`V="7056 \mathcode`W="7057 \mathcode`X="7058
\mathcode`Y="7059 \mathcode`Z="705A

\def\spacedmath#1{\def\packedmath##1${\bgroup\mathsurround=0pt ##1\egroup$}%
\mathsurround#1 \everymath={\packedmath}\everydisplay={\mathsurround=0pt }}

\def\nospacedmath{\mathsurround=0pt \everymath={}\everydisplay={} }

\def\eqalign#1{\null\,\vcenter{\openup\jot\m@th\ialign{
\strut\hfil$\displaystyle{##}$&$\displaystyle{{}##}$\hfil
&&\quad\strut\hfil$\displaystyle{##}$&$\displaystyle{{}##}$\hfil
\crcr#1\crcr}}\,}

\def\system#1{\left\{\null\,\vcenter{\openup1\jot\m@th
\ialign{\strut$##$&\hfil$##$&$##$\hfil&&
        \enskip$##$\enskip&\hfil$##$&$##$\hfil\crcr#1\crcr}}\right.}


\let\@ldmessage=\message

\def\message#1{{\def\pc{\string\pc\space}%
                \def\'{\string'}\def\`{\string`}%
                \def\^{\string^}\def\"{\string"}%
                \@ldmessage{#1}}}



\def\up#1{\raise 1ex\hbox{\smallf@nt#1}}


\def\cf{{\it cf}} 

\def\qed{\raise -2pt\hbox{\vrule\vbox to 10pt{\hrule width 4pt
                 \vfill\hrule}\vrule}}

\def\cqfd{\unskip\penalty 500\quad\vrule height 4pt depth 0pt width
4pt\medbreak}

\def\virg{\raise .4ex\hbox{,}}   


\def\build#1_#2^#3{\mathrel{
\mathop{\kern 0pt#1}\limits_{#2}^{#3}}}


\def\boxit#1#2{%
\setbox1=\hbox{\kern#1{#2}\kern#1}%
\dimen1=\ht1 \advance\dimen1 by #1 \dimen2=\dp1 \advance\dimen2 by #1
\setbox1=\hbox{\vrule height\dimen1 depth\dimen2\box1\vrule}%
\setbox1=\vbox{\hrule\box1\hrule}%
\advance\dimen1 by .4pt \ht1=\dimen1
\advance\dimen2 by .4pt \dp1=\dimen2  \box1\relax}

\def\date{\the\day\ \ifcase\month\or janvier\or f\'evrier\or mars\or avril\or
mai\or juin\or juillet\or ao\^ut\or septembre\or octobre\or novembre\or
d\'ecembre\fi \ {\old \the\year}}

\def\dateam{\ifcase\month\or January\or February\or March\or April\or May\or
June\or July\or August\or September\or October\or November\or December\fi \
\the\day ,\ \the\year}

\def\crog{{\vrule height 2.57mm depth 0.85mm width 0.3mm}\kern -0.36mm [}

\def\crod{]\kern -0.4mm{\vrule height 2.57mm depth 0.85mm width 0.3 mm}}

\def\rond{\kern 1pt{\scriptstyle\circ}\kern 1pt}

\def\hfl#1#2{\nospacedmath\smash{\mathop{\hbox to
12mm{\rightarrowfill}}\limits^{\scriptstyle#1}_{\scriptstyle#2}}}

\def\ghfl#1#2{\nospacedmath\smash{\mathop{\hbox to
25mm{\rightarrowfill}}\limits^{\scriptstyle#1}_{\scriptstyle#2}}}

\def\phfl#1#2{\nospacedmath\smash{\mathop{\hbox to
8mm{\rightarrowfill}}\limits^{\scriptstyle#1}_{\scriptstyle#2}}}

\def\va{vari\'et\'e ab\'elienne}

\def\vapp{vari\'et\'e ab\'elienne principalement polaris\'ee}

\def\pa{\S\kern.15em}

\def\Z{\hbox{\bf Z}}

\def\Q{\hbox{\bf Q}}
\def\C{\hbox{\bf C}}

\def\cf{{\it cf.\/}}

\def\Ext{\mathop{\rm Ext}\nolimits}

\def\isom{\simeq}

\def\loc{{\it loc.cit.\/}}
\def\long{\mathop{\rm long}\nolimits}
\def\lra{\longrightarrow}
\def\llra{\nospacedmath\hbox to 10mm{\rightarrowfill}}
\def\lllra{\nospacedmath\hbox to 15mm{\rightarrowfill}}

\def\NS{\mathop{\rm NS}\nolimits}

\def\ra{\rightarrow}

\def\Sing{\mathop{\rm Sing}\nolimits}

\def\theo{th\'eor\`eme}

\def\tv{\tvi\vrule}
\def\tvi{\vrule height 12pt depth 5pt width 0pt}
\def\tx{\kern -1.5pt -}

\def\JT{(JC,\theta)}

\def\cc#1{\hfill\kern .7em#1\kern .7em\hfill}

\def\dra{\ra\kern -3mm\ra}
\def\ldra{\lra\kern -3mm\ra}

\def\og{\leavevmode\raise.3ex\hbox{$\scriptscriptstyle\langle\!\langle$}}
\def\fg{\leavevmode\raise.3ex\hbox{$\scriptscriptstyle\,\rangle\!\rangle$}}

\catcode`\@=12

\showboxbreadth=-1  \showboxdepth=-1


\baselineskip=12pt
\spacedmath{1.3pt}
\font\eightrm=cmr10 at 8pt
\font\gros=cmbx10 at 12pt
\font\bfit=cmbxti10 at 8pt
 \overfullrule=0mm
\def\ppav{principally polarized abelian variety}

\def\ssv{sous-vari\'et\'e}
\def\XT{(X,\theta)}
\parskip=1mm
\def\ind{\hskip 5mm\relax}

G\'eom\'etrie alg\'ebrique/{\it Algebraic Geometry}
\smallskip \ctitre  {\gros  Sous-vari\'et\'es de codimension $2$

 d'une \va }
\endctitre\vskip 1.6mm

\centerline{Olivier {\pc DEBARRE}}
\footnote {}{\eightpoint\ind Note pr\'esent\'ee par Jean-Pierre {\pc SERRE}}
\bigskip
\baselineskip=9pt {\parindent=5mm\narrower{\eightpoint  {\bfit R\'esum\'e} --
{\it Soit $\XT$ une  vari\'et\'e ab\'elienne principalement polaris\'ee
complexe. On \'etudie les sous-vari\'et\'es $V$ de $X$ dont la classe est
$d\theta^2/2$ ($d$ entier). Supposons $\XT$ g\'en\'erale de dimension $\ge 11$;
lorsque $d$ est impair, $V$ est (tr\`es) singuli\`ere; lorsque $d$ est pair
$\le 16$ et que $V$ est lisse, $V$ est (presque) une intersection
compl\`ete.}}\par}
\bigskip

{\eightpoint\centerline{\bf Subvarieties of codimension $2$ of an abelian
variety}}
\bigskip

{\parindent=5mm\narrower{\eightpoint  {\bfit Abstract} -- {\it Let $\XT$ be a
complex \ppav . Subvarieties of $X$ with class $d\theta^2/2$ ($d$ integer) are
studied. Assume $\XT$ is general of dimension $\ge 11$; when
$d$ is odd, $V$ is (very) singular; when $d$ is even $\le 16$ and $V$ is
smooth,
$V$ is (almost) a complete intersection.}}\par}

\vskip 8mm\parskip=1mm\baselineskip=13pt
 \ind Soit $\XT$ une vari\'et\'e ab\'elienne principalement polaris\'ee
 complexe de dimension $n$. Pour tout entier  $j$, on note $\theta_j$ la
classe de cohomologie enti\`ere $\theta^j/j!$; on s'int\'eresse aux \ssv s de
$X$ de classe
$d\theta_2$ ($d$ entier). Quand $\XT$ est {\it g\'en\'erale}, toute
sous-vari\'et\'e de
$X$ de codimension
$2$ a cette propri\'et\'e ([8]). On trouve de telles \ssv s avec
 $d$ {\it impair} sur les jacobiennes de courbes; par ailleurs, si $X$
contient une courbe irr\'eductible $C$ de classe
$c\theta_{n-1}$ avec
$c$ impair (il existe des exemples qui ne sont pas des jacobiennes), la somme
de $n-2$ copies de $C$ dans $X$ a pour classe $c^{n-2}\theta_2$ et son lieu
singulier est de dimension $\ge n-6$ ([3], prop. 4.7). Pour
$\XT$ g\'en\'erale, la question de l'existence pour $d$ impair peut se
reformuler de fa\c con plus frappante: {\it $\theta_2$ est-elle la classe d'un
cycle alg\'ebrique}? Sans r\'epondre \`a cette question, on montre dans cette
Note que, {\it pour $\XT$ g\'en\'erale de dimension
 $n\ge 7$,  toute \ssv\ de $X$ de classe $d\theta_2$, $d$ impair, est
singuli\`ere et son lieu singulier est de dimension $\ge n/3-2\log_2n$}. On
montre aussi que {\it sur une \vapp\ $\XT$ quelconque de dimension $\ge 8$, il
n'y a pas de fibr\'e vectoriel topologique de rang $2$ et seconde classe de
Chern
$d\theta_2$, $d$ impair}.

\ind Pour $d$ {\it pair}, on cherche un analogue de la conjecture de
Hartshorne sur les \ssv s de l'espace projectif. On dira qu'une \ssv\ de $X$
de codimension $2$ est une {\it intersection presque compl\`ete} si c'est le
lieu des z\'eros d'une section d'un fibr\'e sur $X$ extension de deux fibr\'es
en droites. On montre que, {\it pour $\XT$ g\'en\'erale de dimension $\ge
11$,  toute \ssv\ lisse de $X$ de classe $e\theta^2$, $e\le 8$, est une
intersection presque compl\`ete}.
\smallskip
\ind 1. {\pc PR\'ELIMINAIRES} -- Soit $A^{\bullet}$ une alg\`ebre gradu\'ee
sur $\Q$. Si $c_1\in A^1$ et
$c_2\in A^2$, on note $p_n(c_1,c_2)$ le terme de degr\'e $n$ du caract\`ere de
Chern: si on \'ecrit formellement $c_1=\lambda+\mu$ et $c_2=\lambda\mu$,  on a
$p_n( c_1,c_2)=(\lambda^n+\mu^n)/n!$. Fixons $\theta\in A^1$; pour tous
entiers $j$, $d$ et $t$, et  $\ell\in A^1$, on pose $\ell_j=\ell^j/j!$  et
$p_n(\ell,d,t)=p_n(\ell-2t\theta,d\theta_2-t\ell\theta+t^2\theta^2)$. Pour
tout r\'eel $x$, on note $[x]$ le plus grand entier $\le x$ et
 $\lceil x\rceil$le plus petit entier $\ge x$. Soient $R$ le localis\'e de
$\Z$ en l'id\'eal premier
 $(2)$, vu comme un sous-anneau de $\Q$, et $2R$ son id\'eal maximal. Un
calcul fastideux montre la formule suivante
$$p_n(\ell,d,t)= \sum_{0\le h\le n\atop 0\le j\le h/2} (-1)^{h-j} t^{h-2j}
2^{-\lceil (n-h)/2\rceil -j+1} d^j \ell_{n-h}\theta_h {h\choose 2j}
\rho_{j,n-h}\ ,\leqno (1)$$ o\`u $\rho_{j,n-h}$ est dans $2R$, sauf lorsque
$h=n$, auquel cas il vaut $1$.
\smallskip
\ind Soient  $\XT$ une \vapp\ de dimension $n$,  $\Theta$ un diviseur th\^eta
et $X'$ une intersection compl\`ete g\'en\'erale dans $X$ de $s$ \'el\'ements
de $|3\Theta |$.
\smallskip {\pc PROPOSITION} 2.-- {\it Soit $\ell\in H^2(X,\Z )$. S'il existe
un fibr\'e vectoriel topologique $E$ de rang $2$ sur $X'$ avec
$c_1(E)=\ell_{|X'}$ et $c_2(E)=(d\theta_2)_{|X'}$, le rationnel $\Delta^s
p_n(\ell,d)=\sum_{i=0}^s (-1)^i{s\choose i} p_n(\ell,d,3i)$ est entier.}

\ind Le caract\`ere de Chern $ch(E)$ est la restriction \`a $X'$ d'une classe
$ch_X(E)\in H^{\bullet}(X,\Z)$. L'analogue du
\theo\ de Riemann-Roch pour les fibr\'es topologiques ([6], th.~24.5.4) montre
que $\int_{X'}td(TX')ch(E)$ est entier; or ce nombre vaut
$\sum_{i=0}^s(-1)^s{s\choose i}\int_X ch_X(E)ch\bigl({\cal
O}_X(-3i\Theta)\bigr)$. Puisque $\int_X ch_X(E)ch\bigl({\cal
O}_X(-3i\Theta)\bigr)=p_n(\ell,d,3i)$, ceci d\'emontre le lemme.\cqfd

\smallskip
\ind Soient maintenant $V$ une \ssv\ de $X$ de classe
$d\theta_2$ et $s$ un entier
$>\dim\bigl(\Sing (V)\bigr)$, de sorte que l'intersection $V'=V\cap X'$ est
lisse, contenue dans le lieu lisse de $V$.
\smallskip {\pc PROPOSITION} 3.-- {\it Si $\NS (X)\isom\Z\theta$ et $0\le s\le
n-6$, il existe un fibr\'e vectoriel (alg\'ebrique) $E$ de rang $2$ sur $X'$
v\'erifiant
$c_2(E)=(d\theta_2)_{|X'}$.}

\ind Puisque $X$ est simple, le fibr\'e normal $N_{V'/X'}$, restriction de
 $N_{V/X}$ \`a $V'$, est
 ample (\cf\ [3], prop. 1.1). Puisque
$s\le n-6$, le \theo\ 4.5 de \loc\ entra\^ine qu'il existe un fibr\'e en
droites
$L$ sur $X$ dont la restriction \`a $V'$ est $\wedge^2 N_{V'/X'}$.
L'hypoth\`ese
$\NS (X)\isom\Z\theta$ assure que $H^2(X',L^*_{|X'})=0$. On peut donc
effectuer la construction de Serre ([4], p. 153).\cqfd
\smallskip
 {\it Remarque} 4. -- Dans une jacobienne g\'en\'erale
$\JT$ de dimension $n\ge 6$, la \ssv\ $V=W_{n-2}(C)$ a pour classe $\theta_2$
et son lieu singulier est de dimension $n-6$. Pour
$3\le n\le 5$, $V$ est lisse, mais n'est pas sous-canonique dans $X$: on a
$c_1(\omega_V)=x+\theta$, o\`u
$x$ est la classe de $W_{n-3}(C)$ dans $V$ ([7]).
\medskip

\ind 2. {\pc SOUS}-{\pc VARI\'ET\'ES} {\pc DE} {\pc CLASSE}  $d\theta_2$, $d$
{\pc IMPAIR} -- Soient comme ci-dessus $(X,\theta)$ une \vapp\ de dimension
$n$,
$\Theta$ un diviseur th\^eta et $X'$ une intersection compl\`ete g\'en\'erale
dans
$X$ de $s$ \'el\'ements de $|3\Theta |$.
\smallskip
 {\pc TH\'EOR\`EME} 5.-- {\it Pour $n\ge 4$, $n\ne 6,7$ et $s\le  n-4[{n\over
4}]$, il n'existe pas sur  $X'$ de fibr\'e vectoriel topologique de seconde
classe de Chern
$(d\theta_2)_{|X'}$, $d$ impair.}

\ind Soit $E$ un tel fibr\'e; d'apr\`es le \theo\ de Lefschetz, il existe
$\ell\in H^2(X,\Z )$ avec $c_1(E)=\ell_{|X'}$. Par la proposition~2,
$\Delta^{s'} p_n(\ell,d)$ est entier pour $s\le s'\le n$, en particulier pour
$s'= n-4[{n\over 4}]$. En utilisant la formule (1), on montre qu'il existe
 $\nu(n)>0$ tel que l'entier (pair) $2^{\nu(n)}\Delta^{s'}p_n(\ell,d)$ soit
congru \`a $d$ modulo $2R$, ce qui contredit le fait que $d$ est impair.\cqfd

\smallskip {\pc TH\'EOR\`EME} 6.-- {\it Soit $(X,\theta)$ une
\vapp\ de dimension $n\ge 7$ v\'erifiant $\NS (X)\isom \Z\theta$. Toute
\ssv\ $V$ de $X$ de classe $d\theta_2$, $d$ impair, est singuli\`ere et son
lieu singulier est de dimension
$\ge n/3-2\log_2n$.}

\ind Supposons $n\ge 30$; posons $t_n=\lceil n/3-2\log_2n\rceil$ et supposons
$\dim\bigl(\Sing (V)\bigr)<t_n$. Posons $m=[(n-t_n)/4]$ et $s=n-4m$; on a
$s>\dim\bigl(\Sing (V)\bigr)$ et $s\le n-6$. Soit $E$ le fibr\'e vectoriel sur
$X'$ fourni par la proposition~3. On a
$c_1(E)=a\theta_{|X'}$, et, par la proposition~2,  $\Delta^s p_n(a\theta,d)$
est  entier. On peut v\'erifier que $s\le 2m-\log_2(6m-4)$, ce qui contredit
le lemme suivant.
 \smallskip {\pc LEMME} 7.-- {\it Si $m>0$ et $s\le 2m-\log_2(6m-4)$,
$\Delta^s p_{4m+s}(a\theta,d)$ n'est pas entier.}

\ind La formule (1) exprime $\Delta^s p_{4m+s}(a\theta,d)$ comme une somme
dont une analyse minutieuse des termes montre qu'il n'y en a qu'un de
valuation en $2$ minimale: il correspond
\`a $h=4m+s$ et $j=2m$ et est le produit de $((4m+s)!/ (4m)!)\ 2^{-2m+1}$ par
une unit\'e de $R$. L'hypoth\`ese sur $s$ entra\^ine que
$2^{2m-1}$ ne divise pas
$(4m+s)!/(4m)!$, de sorte que cette valuation est strictement n\'egative.\cqfd

\ind Pour les petites valeurs de $n$, des arguments similaires donnent les
bornes inf\'erieures suivantes $s_n$ sur la dimension du lieu singulier de $V$:

{\eightpoint
$$\vbox{\offinterlineskip \halign{\tv\hskip.5mm \tv#&\cc{$#$}&  \tv#&
\cc{$#$}&  \tv\hskip.5mm \tv#&\cc{$#$}&  \tv#&\cc{$#$}& \tv\hskip.5mm
\tv#&\cc{$#$}& \tv#&\cc{$#$}& \tv\hskip.5mm  \tv#&\cc{$#$}&
\tv#&\cc{$#$}&  \tv\hskip.5mm \tv#&\cc{$#$}&  \tv#&\cc{$#$}& \tv\hskip.5mm
\tv#&\cc{$#$}&  \tv#&\cc{$#$}&  \tv\hskip.5mm \tv#&\cc{$#$}& \tv#&
\cc{$#$}& \tv\hskip.5mm \tv#&\cc{$#$}& \tv#&
\cc{$#$}& \tv\hskip.5mm \tv#\cr
 \noalign{\hrule}\tvi
 &n&& s_n  &&n&&  s_n &&n&& s_n  &&n&& s_n  &&n&& s_n &&n&& s_n &&n&& s_n
&&n&& s_n &\cr
 \noalign{\hrule}\tvi
 &7&& 0 &&  10&& 2 &&  13&& 1 &&  16&& 2 &&  19&& 5 &&  22&& 6 &&  25&& 5&&
28&& 8&\cr
 &8&& 0 &&  11&& 2 &&  14&& 2 &&  17&& 1 &&  20&& 4 &&  23&& 7 &&  26&& 6&&
29&& 9&\cr
 &9&& 1 &&  12&& 2 &&  15&& 4 &&  18&& 2 &&  21&& 5 &&  24&& 6 &&  27&& 9&&
30&& 10&\cr
\noalign{\hrule} }}$$}

\ind En particulier, $V$ est singuli\`ere pour $n\le 30$. Puisque $t_n\le 0$
pour $n<30$ et
$t_n\ge 0$ pour $n\ge 30$,  cela montre le \theo .\cqfd

\smallskip
\ind 3. {\pc SOUS}-{\pc VARI\'ET\'ES} {\pc DE} {\pc CLASSE}  $d\theta_2$, $d$
{\pc PAIR} -- Soit $\XT$ une \vapp ; une \ssv\ $V$ de $X$ de codimension $2$
est une {\it intersection presque compl\`ete} si c'est le lieu des z\'eros
d'une section
$\sigma$ d'un fibr\'e
$E$ extension de fibr\'es en droites $M'$ et $M$ sur $X$. Si
$M\not\isom M'$, l'extension est scind\'ee et $V$ est une intersection
compl\`ete. Si $M\isom M'$, que $e_E\in H^1(X,{\cal O}_X)\isom
\Ext^1(M,M)$ est la classe de l'extension, et que $s$ est l'image de
$\sigma$ dans $H^0(X,M)$,
$V$ est contenue dans le diviseur $Z(s)$ de
$s$ et est le lieu des z\'eros d'une section de $M_{|Z(s)}$. Une telle  section
peut s'\'ecrire
$s'+Ds$, avec $s'\in H^0(X,M)$ et $D\in H^0(X,TX)$ (\cf\ [1] pour la
d\'efinition de $Ds$), et
$V$ peut \^etre d\'ecrite par les \og\'equations\fg\ $s=s'+Ds=0$, o\`u $D\smile
c_1(M)\in\C e_E$.
\smallskip
 {\pc TH\'EOR\`EME} 8.-- {\it Soient $(X,\theta)$ une \vapp\ v\'erifiant $\NS
(X)\isom\Z\theta$, et $V$ une \ssv\ lisse de
 $X$ de classe $e\theta^2$, avec $ e\le 8$. On suppose la dimension $n$ de $X$
minor\'ee de la fa\c con suivante

{\eightpoint
$$\vbox{\offinterlineskip \halign{\tv#&\ \cc{$#$}\ &\tv\hskip.5mm
\tv#&\cc{$#$}&\tv#&\cc{$#$}&\tv#&\cc{$#$}&\tv#&\cc{$#$}&\tv#
&\cc{$#$}&\tv#&\cc{$#$}&\tv#&\cc{$#$}&\tv#&\cc{$#$}&\tv#&\cc{$#$}&\tv#&\cc{$#$}
&\tv#&\cc{$#$}&\tv#&\cc{$#$}&\tv#&\cc{$#$}&\tv#&\cc{$#$}&\tv#&\cc{$#$}&\tv#\cr
\noalign{\hrule} &e&&1&&2&&3&&4&&5&&6&&7&&8&\cr
\noalign{\hrule} &n\ge&&6&&6&&7&&6&&8&&7&&11&&8&\cr
\noalign{\hrule} }}$$}
\ind Si $e\ge 4$, on suppose de plus les diviseurs th\^eta non singuliers en
codimension $3$; alors $V$ est une intersection presque compl\`ete.}

\ind Comme $X$ est simple, le fibr\'e normal $N$ de $V$ dans $X$ est ample;
puisque $n\ge 6$, il existe d'apr\`es [3], th. 4.5
 un entier $a$ et un diviseur th\^eta $\Theta'$ tels que $\omega_V\isom{\cal
O}_V(a\Theta')$. La construction de Serre fournit un fibr\'e vectoriel
$E$ sur $X$ de classes de Chern $a\theta$ et $e\theta^2$, dont la restriction
\`a $V$ est $N$. Comme
 $N$ est ample, on a  $a^2> 4e\cos^2{\pi\over n-1}$ ([10]) et, vu les valeurs
de $n$ et $e$, $a^2\ge 4e$.
\smallskip {\pc LEMME} 9.-- {\it Il existe un entier $b\ge a/2$ v\'erifiant
$b(a-b)\le e$, et un diviseur th\^eta
$\Theta$, tels que
$H^0\bigl( X,E(-b\Theta)\bigr)\ne 0$. Si $b(a-b)= e$, $V$ est une intersection
presque compl\`ete.}

\ind Si $a^2>4e$, le fibr\'e $E$ n'est pas $\theta$\tx semi-stable ([2], [9],
th.~4.3). Si $a^2=4e$, $a$ est pair et les classes de Chern de
$E(-(a/2)\Theta')$ sont nulles. Comme $\pi_1(X)$ n'a pas de repr\'esentation
irr\'eductible de dimension
$2$, le fibr\'e $E(-(a/2)\Theta')$ n'est pas
$\theta$\tx stable ([9], th.~5.1). Dans les deux cas, il existe un entier $b\ge
a/2$, un diviseur th\^eta
$\Theta$ et une section non nulle $s$ de
$E(-b\Theta)$. Si $b$ est maximal avec cette propri\'et\'e, $Z(s)$ est vide ou
de codimension
$2$. On en d\'eduit
$0\le c_2\bigl( E(-b\Theta)\bigr)\cdot\theta^{n-2}=\bigl( e -  b(a-b)\bigr)
n!$; s'il y a
\'egalit\'e, $Z(s)$ est vide et $E(-b\Theta)$ est extension de deux fibr\'es
en droites, d'o\`u le lemme.\cqfd
\smallskip
\ind Supposons $b(a-b)<e$; si $a-b\ge 2$, on a $a^2/4\ge e>2b\ge 2\lceil
a/2\rceil$, d'o\`u $a\ge 6$ et $e>2(a-2)\ge 8$, ce qui contredit
l'hypoth\`ese. D'autre part, la suite exacte de la construction de Serre (\cf\
[4]) donne
 $H^0\bigl( X,{\cal I}_V ( a\Theta'-b\Theta)  \bigr)\ne 0$. On en d\'eduit
$b=a-1$, de sorte que $a-1<e\le a^2/4$ et $e\ge a\ge 4$, et qu'un diviseur
th\^eta
$\Theta''$ contient
$V$. Il est non singulier en codimension $3$, donc localement factoriel ([5]),
et  $V$ est un diviseur de Cartier dans $\Theta''$. Comme la classe de
$V$ est $e\theta^2$, on a $\omega_V\sim (e+1)\theta_{|V}$, donc
$a=e+1$, contradiction. On a donc $e=b(a-b)$, et le lemme 9 entra\^ine le
\theo .\cqfd
\smallskip
\ind{\eightrm  Financ\'e en partie par N.S.F. Grant DMS 94-00636 et le projet
europ\'een HCM
\og Algebraic Geometry in Europe\fg (AGE), contrat CHRXCT-940557.
\smallskip Note remise et accept\'ee le 18 septembre 1995}\vskip 5mm
\parindent=2mm {\eightpoint{\pc R\'EF\'ERENCES} {\pc BIBLIOGRAPHIQUES}
\smallskip
\baselineskip=9pt

 [1] A. {\pc BEAUVILLE} et O. {\pc DEBARRE}, {\it Une relation entre deux
approches du probl\`eme de Schottky}, Invent. Math. {\bf 86} (1986), 195--207.

 [2] F. A. {\pc BOGOMOLOV}, {\it Holomorphic tensors and vector bundles on
projective varieties}, Math. USSR Izvestija {\bf 13} (1979), 499--555.

 [3] O. {\pc DEBARRE}, {\it Fulton-Hansen and Barth-Lefschetz Theorems for
Subvarieties of Abelian Varieties}, J. f\"ur die reine und andgewandte
Mathematik (1995), \`a para\^itre.

 [4]	D. {\pc FERRAND},  Construction de fibr\'es de rang deux,  in {\it Les
\'equations de Yang-Mills}, S\'eminaire E.N.S. 1977-1978, Ast\'erisque 71-72,
1980.

 [5] A. {\pc GROTHENDIECK}, Cohomologie locale des faisceaux coh\'erents et
\theo s de Lefschetz locaux et globaux (SGA 2), Masson et North Holland, Paris
Amsterdam, 1968.

 [6] F. {\pc HIRZEBRUCH},  Topological methods in algebraic geometry,
Springer Verlag, Berlin Heidelberg New-York, 3\up{e} \'ed., 1966.

 [7] I. {\pc MACDONALD}, {\it Symmetric Products of an Algebraic Curve},
Topology {\bf 1} (1962), 319--343.

 [8] T. {\pc MATTUCK}, {\it Cycles on abelian varieties},  Proc. Amer. Math.
Soc. {\bf 9} (1958), 88--98.

 [9]	V. B. {\pc MEHTA} et A. {\pc  RAMANATHAN}, {\it Restriction of stable
sheaves and representation of the fundamental group}, Invent. Math. {\bf 77}
(1984), 163--172.

 [10]	M. {\pc SCHNEIDER},  Submanifolds of projective space with semistable
normal bundle,  in {\it Several Complex Variables, Proc. of the 1981 Hangzhou
Conf.}, J.J.Kohn, Q.-K. Lu, R. Remmert et Y.-T. Siu \'ed., Birkh\"auser,
Boston 1984, 151--160.

\smallskip
\hfill\hbox to 2cm{\hrulefill}\parskip=0cm

\hfill{\it Math\'ematique, Universit\'e Louis Pasteur, 7 rue Ren\'e Descartes,
67084 Strasbourg C\'edex, France.

\hfill adresse \'electronique:} debarre@math.u-strasbg.fr
\par}
\bye